\documentclass[12pt]{amsart}
\usepackage{amsmath, amssymb}
\usepackage{amsmath,amsfonts,amssymb,amsthm,amscd} 
\usepackage{enumerate}
\newcommand{\inn}{\operatorname{Inn}}

\newcommand{\sym}{\operatorname{Sym}}

\newcommand{\syl}{\operatorname{Syl}}

\newcommand{\stab}{\operatorname{Stab}}

\newcommand{\SL}{\mathrm{SL} }
\newcommand{\PSL}{\mathrm{PSL}}

\newcommand{\G}{\mathrm{G}}

\newcommand{\B}{\mathrm{B}}

\newcommand{\GL}{\mathrm{GL}}

\newcommand{\A}{\AA}

\newcommand{\M}{\mathrm{M}}

\newcommand{\Q}{\mathrm{Q}}

\renewcommand{\AA}{{\mathcal  A}}
\def \Aut {{\mathrm{Aut}} }

\def \AGL {\mathrm{AGL}}
\def \<{\langle }
\def \>{\rangle }

\def \inv {^{-1}}

\def \GL {\mathrm{GL} }
\def \sym {\mathrm{Sym} }
\def \syl {\mathrm{Syl} }

\renewcommand{\leq}{\leqslant}
\renewcommand{\geq}{\geqslant}
\renewcommand{\le}{\leqslant}

\def \b{\beta }
\def \a{\alpha }
\def \g{\gamma }

\def \G{\Gamma }

\input cyracc.def

\newtheorem*{thmA}{Theorem A}
\newtheorem*{thmB}{Theorem B}

\newtheorem{lemma}{Lemma}
\newtheorem{thm}[lemma]{Theorem}

\title{A $3$-Local Characterization of $\mathrm{M}_{12}$ and $\SL_3(3)$}
\author{Sarah Astill} \address{School of Mathematics\\
University of Birmingham\\
Edgbaston\\
Birmingham B15 2TT\\
United Kingdom} \email{astills@maths.bham.ac.uk}

\author{Chris Parker}
\address{School of Mathematics\\
University of Birmingham\\
Edgbaston\\
Birmingham B15 2TT\\
United Kingdom} \email{c.w.parker@bham.ac.uk}

\begin{document}\maketitle
\begin{abstract}
We identify the sporadic simple group $\mathrm{M}_{12}$ and the
simple group $\mathrm{SL}_3(3)$ from some part of their $3$-local
structure and give a graph theoretic analogue of the resulting
theorem.
\end{abstract}

Let  $G_1=\mathrm{SL}_3(3)$ and set $A_1
=\mathrm{Stab}_{G_1}(\<(1,0,0)\>)$, $B_1
=\mathrm{Stab}_{G_1}(\<(1,0,0),(0,1,0)\>)$ and $C_1 = A_1\cap B_1$.
Then $A_1 \cong B_1 \cong \AGL_2(3)$ and $C_1= N_{G_1}(Z(S_1))$
where $S_1 \in \syl_3(C_1)$. Now consider $G_2=\mathrm{M}_{12}$, the
Mathieu group of degree 12, acting on the Steiner system $\mathcal
S= S(12,6,5)$ of $\Omega=\{1,\hdots,12\}$. A linked three is a
partition of $\Omega$ into four subsets of size three  such that the
union of any two of the subsets is a hexad in $\mathcal S$. Assume
that $\mathcal T$ is a linked three and $\mathbf t$ is a subset from
$\mathcal T$.
 Let $A_2=\mathrm{Stab}_{G_2}(\mathbf t)$, $B_2=\mathrm{Stab}_{G_2}(\mathcal T)$ and  $C_2 = A_2\cap
B_2$. Then  $A_2 \cong B_2 \cong \AGL_2(3)$ and $C_2= N_{G_2}(Z(S_2))$ where $S_2 \in \syl_3(C_2)$.
Thus from  the perspective of the normalizers of $3$-subgroups the simple group $\SL_3(3)$ and the
sporadic simple group $\mathrm{M}_{12}$ are indistinguishable. Our main theorem is as follows.

\begin{thmA}\label{maintheorem}
Let $G$ be a finite group, $A, B \leq G$, $A \neq B$ and $S \in \mathrm{Syl}_3(A\cap B)$. Assume
that $G = \langle A,B\rangle$ and $A \cong B\cong \mathrm{AGL}_2(3)$. If $N_G(Z(S))\leq A$, then
$G\cong\mathrm{SL}_3(3)$ or $G \cong \mathrm{M}_{12}$.
\end{thmA}

Our proof of Theorem~A is relatively elementary. It uses a mixture
of local group theoretic methods, character theory in the guise of a
theorem of Feit and Thompson and a coset enumeration package. Of
course our proof does not use the classification of finite simple
groups. We expect that the result will be useful in the ongoing
project to understand the groups of local characteristic $p$
\cite{MSS-overview}. As is often the case in $p$-local
characterizations for odd primes $p$, our proof is centered around
the determination of the structure of an involution centralizer. Let
$G$, $A$, $B$, $C$ and $S$ be as in the statement of the theorem and
let $t \in A$ be an involution. Then it turns out that $C_G(t)$ has
Sylow $3$-subgroups of order three. This means that
 local group theoretic methods cannot be used to restrict the structure of $C_G(t)$.
However, we can show that the Sylow $3$-subgroup of $C_G(t)/\<t\>$
is self-centralizing. The structure of finite groups with a
self-centralizing subgroup of order three have been described by
Feit and Thompson in \cite{FeitThompson}. In fact a generalization
of the Feit--Thompson Theorem by Mazurov \cite{Mazurov} would allow
us to strengthen the theorem to allow for, not necessarily finite,
groups $G$ for which $\<Z(S),Z(S)^g\>$ is finite for any $g \in G$.
 A lemma originally due to Burnside  (see Lemma~\ref{commutator relation})
provides us with a list of possible group relations satisfied by
certain commutators of elements in $C_G(t)$. By considering amalgams
and using the Goldschmidt Lemma we give a presentation for an
infinite group $F$ which has any  group satisfying the hypothesis of
Theorem~A as a quotient. We then create quotient groups of $F$
defined by the possible group relations given  by the Burnside
Lemma. Somewhat miraculously just two of these groups are
non-trivial and each has finite order. Moreover the group orders
coincide with the orders of our target groups. At this stage we
conclude that Theorem~A is true.

We emphasize that the   Feit--Thompson Theorem only applies  when we
have a self-centralizing subgroup of order three. Thus the
techniques used here do not apply to a similar set-up where $A \cong
B \cong \AGL_2(p)$ for odd primes $p>3$. On the other hand, in
\cite{BN-pairs} groups $G$ generated by subgroups  $A \cong B \cong
\AGL_2(p^a)$ and $a \ge 2$ and satisfying $N_G(Z(S)) =A\cap B$ where
$S \in \syl_3(A\cap B)$ are shown to be isomorphic to $\PSL_3(p^a)$.

In the final section of the paper we prove  a graph theoretic analogue to our main result. Let
$\Gamma_1$ be the graph  whose vertices are the thirteen points together with the thirteen lines in
the projective plane of order three with a point adjacent only to those lines containing it.
Suppose that $\mathcal S=S(12,6,5)$ is the Steiner system on $\Omega = \{1,\dots, 12\}$ mentioned
above. Let $\mathcal R$ be the set of subsets of $\Omega$ of size three and $\mathcal L$ be the set
of linked threes determined by $\mathcal S$. Define $\Gamma_2$ to be the graph with vertex set
$\mathcal R\cup \mathcal L$ and edge set $\{\{\{a,b,c,d\},a\} \mid \{a,b,c,d\} \in\mathcal L, a\in
\mathcal R\}$. Obviously, $\PSL_3(3)\cong \SL_3(3)$ acts as a group of automorphisms of $\Gamma_1$
and $\M_{12}$ acts as a group of automorphisms of $\Gamma_2$. In both cases the action is
transitive on the edges and not on vertices.

\begin{thmB}\label{TheoremB}
Suppose $\Gamma$ is a finite, connected graph and $G\leq
\Aut(\Gamma)$ with $G_\a=\mathrm{Stab}_G(\a) \cong
\mathrm{AGL}_2(3)$ for each vertex $\a$. Suppose further that $G$ is
edge transitive, has two orbits on the set of vertices and that for
some vertex $\alpha$ and each non-trivial $z \in O_3(G_\a)$ the
subgraph of $\Gamma$ fixed by $\<z\>$ is a tree which contains at
least one edge. Then $G \cong \mathrm{PSL}_3(3)$ or $G\cong
\mathrm{M}_{12}$ and $\Gamma$ is isomorphic to $\Gamma_1$ or
$\Gamma_2$ respectively.
\end{thmB}

We remark that in principle the coset enumerations in this paper can
be carried out by hand. This would of course make Theorem~A
independent from computer calculations.

Finally we mention that notation for groups is standard as in
 \cite{Aschbacher} and \cite{atlas}. In particular, for a
group $G$ and a set of primes $\pi$ we denote $O_\pi(G)$ to be the
largest normal subgroup of $G$ which is a $\pi$-group.

\noindent {\bf Acknowledgement.}  The first author acknowledges
financial support from EPSRC in the form of a postgraduate
scholarship.

\section{Preliminaries}

In  this section we present Burnside's Lemma and state the
Feit--Thompson Theorem mentioned in the introduction. We also
present an elementary lemma about $\AGL_2(3)$.

\begin{lemma}\label{commutator relation}
Let $G$ be a group and let $a,b,c \in G$. The following identities
hold. \begin{enumerate}[$(i)$]
        \item $[a,bc]=[a,c][a,b]^c$ and $[ab,c]=[a,c]^b[b,c]$.
        \item If $p$ is a prime and if $G$ is a $p$-group, then $[[a,b],c][[b,c],a][[c,a],b]\in [G,G,G,G]$.
      \end{enumerate}
\end{lemma}
\begin{proof}
The   identities in (i) are easy to check. Part (ii) is Lemma 5.6.1 (iv) in \cite[p.
209]{gorenstein}.
\end{proof}

The following lemma is well-known and dates back to Burnside. It can
be found in his book \cite[p. 90-91]{BurnsideBook}. We present a
detailed  proof of the result since we want to record a certain
identity which will be of major importance to us in the proof of
Theorem~A. Our proof is modelled on one that can be found in
\cite[Theorem 8.1]{Higman} which is unpublished. We will use the
notation $a \equiv b ~\mathrm{mod}~ H $ when $Ha=Hb$ where $H$ is a
subgroup of a group $G$ and $a,b \in G$.

\begin{lemma}[Burnside]\label{special2^7}
Suppose that $p$ is a prime  and $Q$ is a $p$-group. Assume that ${\zeta}\in \Aut(Q)$ has order
three and $C_Q({\zeta})=1$. Then $Q$ has class at most two and for all $v,w\in Q$,
\[[v,w^{\zeta}] = [v^{\zeta},w] =[v,w]^{{\zeta}^2} .\]
\end{lemma}
\begin{proof}
Let $Q=Q_1 \geq Q_2 \geq \hdots \geq Q_n=1$ be the lower central series of $Q$. Let $i,j\geq 1$ and
let $x \in Q_i$ and $y \in Q_j$. Then $xx^{\zeta}x^{{\zeta}^2}$ is centralized modulo $Q_{i+1}$ so
\[xx^{\zeta}x^{{\zeta}^2} \equiv 1 ~\mathrm{mod}~Q_{i+1}\] and
\[ [xx^{\zeta}x^{{\zeta}^2},y]\in [Q_{i+1},Q_j]\leq Q_{i+j+1}.\]  Since
$Q_{i+j}/Q_{i+j+1}$ is central in $Q/Q_{i+j+1}$, we may use the
commutator relations from  Lemma~\ref{commutator relation} (i) to
get $[x,y][x^{\zeta},y][x^{{\zeta}^2},y] \equiv
[xx^{\zeta}x^{{\zeta}^2},y] \equiv 1 ~\mathrm{mod}~ Q_{i+j+1}$ and
so $[x,y]\equiv [y,x^{\zeta}][y,x^{{\zeta}^2}] ~\mathrm{mod}~
Q_{i+j+1}$. Now swapping the roles of $x$ and $y$ gives $[x,y]\equiv
[x,y^{\zeta}][x,y^{{\zeta}^2}] ~\mathrm{mod}~ Q_{i+j+1}$. Therefore
\[[x^{{\zeta}^2},y][y^{{\zeta}^2},x] \equiv
[x,y^{\zeta}][y,x^{\zeta}] ~\mathrm{mod}~ Q_{i+j+1}\] and is thus centralized by $\zeta$.
 Hence we have $[x,y^{\zeta}]\equiv [x^{\zeta},y] ~\mathrm{mod}~
Q_{i+j+1}$. If we replace $x$ by $x^{\zeta^2} $and $y$ by $y^{\zeta^2}$, we get $[x^\zeta,y]\equiv
[x,y]^{\zeta^2} \mathrm{mod}~ Q_{i+j+1}.$ Thus
\begin{equation}{\label{one}} [x,y^{\zeta}]\equiv [x^{\zeta},y] \equiv [x,y]^{\zeta^2}~\mathrm{mod}~
Q_{i+j+1}.\end{equation}
 In particular,  for all $d,e,f\in Q$, we have
\begin{equation}{\label{one+}}[d,e^{\zeta}]\equiv [d^{\zeta},e]\equiv [d,e]^{{\zeta}^2} ~\mathrm{mod}~Q_3\end{equation} and
\begin{equation}{\label{two}} [[d,e],f^{{\zeta}^2}]\equiv [[d,e]^{\zeta},f^{\zeta}] \equiv [[d,e]^{{\zeta}^2},f]\equiv
[[d^{\zeta},e],f] \equiv [[d,e^{\zeta}],f] ~\mathrm{mod}~Q_4.
\end{equation}
Now let $a,b,c \in Q$. Then Lemma~\ref{commutator relation} (ii)
gives
\begin{equation}\label{three}[[a,b],c][[b,c],a][[c,a],b] \equiv 1
~\mathrm{mod}~Q_4\end{equation} and
\begin{equation}\label{four}[[a^{\zeta},b],c][[b,c],a^{\zeta}][[c,a^{\zeta}],b] \equiv 1
~\mathrm{mod}~Q_4.\end{equation} Conjugating (\ref{three}) by ${\zeta}$ and applying (\ref{two})
gives
\[[[a^{\zeta},b],c][[b,c],a^{{\zeta}^2}][[c,a^{\zeta}],b] \equiv 1
~\mathrm{mod}~Q_4\] which together with (\ref{four}) gives
\[[[b,c],a^{{\zeta}^2}]\equiv [[b,c],a^{\zeta}]
~\mathrm{mod}~Q_4.\] Finally, (\ref{one}) gives us that $[[b,c],a^{\zeta}] \equiv
[[b,c],a]^{{\zeta}^2} ~\mathrm{mod}~Q_4$ and so using (\ref{one}) we get
\[[[b,c],a]^{\zeta} \equiv [[b,c],a^{{\zeta}^2}] \equiv
[[b,c],a^{\zeta}] \equiv [[b,c],a]^{{\zeta}^2} ~\mathrm{mod}~Q_4\] is  fixed by $\zeta$ mod $Q_4$.
Therefore,  $[[b,c],a] \in Q_4$ for all $a,b,c \in Q$. Thus $Q_3 \le Q_4 \leq Q_3$ and we conclude
that $Q_3=1$. In particular,  (\ref{one+}) now implies $[v,w^{\zeta}] =[v^{\zeta},w]=
[v,w]^{{\zeta}^2}$ for all $v,w\in Q$ and this concludes the proof of the lemma.
\end{proof}

The following theorem by Feit and Thompson describing groups with a self-centralizing element of
order three can be used to provide a setting in which to apply our Burnside lemma.

\begin{thm}[Feit--Thompson]\label{Feit-Thompson}
Let $G$ be a finite group containing a subgroup, $X$, of order three such that $C_G(X)=X$. Then one
of the following  hold.
\begin{enumerate}[$(i)$]
    \item $G$ contains a nilpotent normal subgroup, $N$, such that $G/N$
    is isomorphic to $\sym(3)$ or $C_3$.
    \item $G$ contains a  normal  2-subgroup, $N$, such that
    $G/N\cong \mathrm{Alt}(5)$.
    \item $G$ is isomorphic to $\mathrm{PSL}_2(7)$.
\end{enumerate}
\end{thm}

Notice that Burnside's Lemma implies that the subgroup $N$ in parts (i) and  (ii) of
Theorem~\ref{Feit-Thompson} is nilpotent of class at most two. In fact, an often used result of
Higman \cite[Theorem 8.2]{Higman}  shows that the subgroup $N$ in part (iii) is abelian.

We now list some elementary facts about groups isomorphic to $\mathrm{AGL}_2(3)$.
\begin{lemma}\label{facts about AGL23}
Let $X \cong \mathrm{AGL}_2(3)$ and $Z=N_X(S)$ where $S\in
\mathrm{Syl}_3(X)$. Then the following hold:
\begin{enumerate}[$(i)$]
  \item $O_3(X)$ is elementary abelian of order $9$, is the unique non-trivial normal $3$-subgroup of $X$,
  is self-centralizing in $X$ and $X/O_3(X) \cong \GL_2(3)$  acts transitively on the non-trivial elements of
  $O_3(X)$;
  \item if $t \in X$ is an involution with $O_3(X)t \in Z(X/O_3(X))$, then $C_X(t)$ is a complement to $O_3(X)$ in $X$ and, in particular,
  $C_X(t) \cong \GL_2(3)$;
  \item $S$ is extraspecial of order $27$ and exponent $3$ and $C_Z(S)=C_X(S)=Z(S)$ is cyclic of order
  $3$;
  \item $Z= N_X(Z(S))$ has index $4$ in $X$ and is a maximal subgroup of
  $X$; and
  \item $\mathrm{Inn}(Z)\cong Z$, $|\mathrm{Aut}(Z)/\mathrm{Inn}(Z)|=2$, $\mathrm{Aut}(X)=\mathrm{Inn}(X)\cong
  X$ and in particular, $N_{\Aut(X)}(Z) \cong  \inn(Z)\cong Z$.
\end{enumerate}
\end{lemma}
\begin{proof}
These can be easily verified. See \cite{AstillMPhilThesis} for their
proofs.
\end{proof}

\section{The Amalgam}\label{amalgamsection}

An {amalgam} (of rank 2) is a quintuple $(X,Y,Z,\phi_1,\phi_2)$ consisting of three groups $X$, $Y$
and $Z$ and two monomorphisms $\phi_1 : Z \rightarrow X$ and $\phi_2 : Z \rightarrow Y$.

Suppose that $\AA_1=(X_1,Y_1,Z_1,\phi_1,\phi_2)$ and $\AA_2=(X_2,Y_2,Z_2,\theta_1,\theta_2)$ are
amalgams. Then $\AA_1$ and $\AA_2$ have the same weak type  provided there exist isomorphisms $\a :
X_1 \rightarrow X_2$, $\b : Y_1 \rightarrow Y_2$ and $\g : Z_1 \rightarrow Z_2$. If, additionally,
it can be arranged that
 $\mathrm{Im}(\phi_1 \a) =\mathrm{Im}( \g \theta_1)$ and  $\mathrm{Im}(\phi_2 \b) =
\mathrm{Im}(\g \theta_2)$, then we say that $\AA_1$ and $\AA_2$ have the same type. Finally, we say
that  $\AA_1$ and $\AA_2$ are isomorphic if and only if $\phi_1 \a = \g \theta_1$ and $\phi_2 \b =
\g
    \theta_2$. A (faithful) {completion} $(H,\psi_1,\psi_2)$ of the amalgam $\AA =(X,Y,Z,\phi_1,\phi_2)$
is a group $H$  and two monomorphisms $\psi_1 : X \rightarrow H$ and
$\psi_2 : Y \rightarrow H$ such that $\phi_1\psi_1 = \phi_2\psi_2$
and $H = \< X\psi_1,Y\psi_2\> $.  A completion $(F,\Psi_1,\Psi_2)$
of $\AA$ is universal if given any other completion
$(H,\psi_1,\psi_2)$ there is a unique homomorphism $\pi : F
\rightarrow H$ such that $\Psi_i \pi = \psi_i$ for $i=1,2$ (notice
that $\pi$ therefore maps the images of $X$, $Y$ and $Z$ in $F$ to
the images of $X$, $Y$ and $Z$ in $H$). If $G$ is a group and if
there exists monomorphisms $\psi_1:X \rightarrow G$ and
$\psi_2:Y\rightarrow G$ such that $(G, \psi_1,\psi_2)$ is a
completion of $\A$, then we shall also say that the group $G$ is a
completion of $\A$. If, in fact, $(G, \psi_1,\psi_2)$  is a
universal completion then we say that $G$ is a universal completion.
The definition of the universal completion implies that any
completion of $\AA$ is isomorphic to a quotient of the universal
completion of $\AA$ and that the universal completion is itself
unique up to isomorphism. It is easy to check that isomorphic
amalgams have isomorphic universal completions. It is well-known
(see \cite[Theorem 1, p.3]{serre}) that the universal completion of
$\mathcal A$ is isomorphic to the quotient of the free product of
$X$ and $Y$ factored by the normal subgroup generated by the set
$\{(z\phi_1)(z\phi_2)^{-1}\mid z \in Z\}$.

The typical way  that amalgams and completions arise is as follows. Let  $H$ be a group with a
tuple of subgroups $(X,Y,Z)$ such that $Z\leq X$, $Z \leq Y$ and $H =\<X,Y\>$. Then taking $\iota_1
: Z \rightarrow X$, $\iota_2 : Z \rightarrow Y$, $\iota_3:X \rightarrow H$ and $\iota_4:Y
\rightarrow H$ to be inclusion maps, we get an amalgam  $\mathcal A=(X,Y,Z,\iota_1 ,\iota_2)$ which
is of weak type $(X,Y,Z)$ as well as a completion $(H,\iota_3,\iota_4)$ of $\A$. Thus in this case
$H$ is a completion of $\A$.

A fundamental result of Goldschmidt provides a method for calculating the number of isomorphism
classes of an amalgam of a given type. We first need a definition. Suppose that $H \le K$. Then
$\Aut(K,H) = N_{\Aut(K)}(H)/C_{\Aut(K)}(H)$ identified as a subgroup of $\Aut(H)$.

\begin{lemma}[Goldschmidt Lemma]\label{goldschmidt} Suppose  that $\AA=(X,Y,Z,\phi_1,\phi_2)$ is an amalgam and  define
$$A_X =\{\phi_1 \a\phi_1^{-1}\mid \a \in \Aut(X,\phi_1(Z))\}$$ and
$$A_Y =\{\phi_2 \b\phi_2^{-1}\mid \b \in \Aut(Y,\phi_2(Z))\}.$$
  Then there is a one to one correspondence between
$(A_X,A_Y)$-double cosets in $\Aut(Z)$ and isomorphism classes of amalgams of the same type as
$\mathcal A$.
\end{lemma}

\begin{proof} See \cite[(2.7)]{goldschmidt}.
\end{proof}

We now introduce two amalgams which are important in the proof of
Theorem~A.  Let $H \cong \mathrm{SL}_3(3)$ and set
$X=\mathrm{Stab}_H{(\<(1,0,0)\>)}$,
$Y=\mathrm{Stab}_H(\<(1,0,0),(0,1,0)\>)$ and  $Z=  X\cap Y$. Notice
that $Z = N_X(S)$ for $S \in \syl_3(Z)\subseteq \syl_3(X)$. Take
$\iota_1 : Z \rightarrow X$, $\iota_2 : Z \rightarrow Y$, $\iota_3 :
X \rightarrow H$ and  $\iota_4 : Y \rightarrow H$ to be inclusions,
and  define
$$\AA=(X,Y,Z,\iota_1 ,\iota_2).$$ Then $(H,\iota_3,\iota_4)$
is a completion of the amalgam $\AA$.

Now take $H' = \AGL_2(3)$ and let $S' \in \syl_3(H')$. Set $Z'=
N_{H'}(S')$ and let $X'=Y'= H'$. Taking  $\rho_1: Z' \rightarrow
X'$, $\rho_2 : Z' \rightarrow Y'$, $\rho_3: X' \rightarrow H'$ and
$\rho_4 : Y' \rightarrow H'$ as inclusions, we have an amalgam
$$\AA'=(X',Y',Z',\rho_1 ,\rho_2)$$ with a completion
$(H',\rho_3,\rho_4)$.

\begin{lemma}\label{A' has non-trivial O3}
The amalgams $\AA$ and $\AA'$ have the same type but are not
isomorphic. If a group $L$ is a  completion of $\AA'$ then
$O_3(L)\neq 1$.
\end{lemma}
\begin{proof}
Obviously $\AA$ and $\AA'$ have the same weak type. Furthermore, as
the images of $Z$ and $Z'$ are uniquely determined up to conjugacy
in $X$ and $Y$ and $X'$ and $Y'$ respectively, $\AA$ and $\AA'$ have
the same type. Suppose that $(L, \psi_1,\psi_2)$ is a completion of
$\AA'$. In $H'$ we have $1\neq Q = O_3(X'\rho_3)= O_3(Y'\rho_4)$.
Set $R=Q\rho_3^{-1}\rho_1^{-1}=Q\rho_4^{-1}\rho_2^{-1}$. Then
$R\rho_1\trianglelefteq X'$ and $R\rho_2\trianglelefteq Y'$.
Therefore $R\rho_1\psi_1\trianglelefteq X'\psi_1$ and
$R\rho_2\psi_2\trianglelefteq Y'\psi_2$. As $L$ is a completion of
$\AA'$, we have $1\neq R\rho_1\psi_1= R\rho_2\psi_2$ is normal in
$L=\<X'\psi_1, Y'\psi_2\>$. Thus every completion of $\AA'$ contains
a non-trivial normal $3$-subgroup. In particular, since $\SL_3(3)$
is a completion of $\AA$, we have $\AA'$ and $\AA$ are not
isomorphic.
\end{proof}

\begin{lemma}\label{uniqueness} There are exactly two isomorphism classes of amalgams
of weak type $(X,Y,Z)$.
\end{lemma}
\begin{proof} First, as subgroups of $X$ and $Y$ which are isomorphic to $Z$ are uniquely determined up to
conjugacy, any two amalgams of weak type $(X,Y,Z)$ in fact have the
same type. We can therefore use Goldschmidt's Lemma  to determine
the number of isomorphism classes of amalgams of weak type
$(X,Y,Z)$. We identify $Z$ with a subgroup $Z^*$ of $X$.  By
Lemma~\ref{facts about AGL23} (v), $N_{\Aut(X)}(Z) \cong
\inn(Z)\cong Z$ and $C_X(Z^*)=1$. Therefore $A_X = \inn(Z)\cong
\inn(Z^*)$ and similarly, $A_Y =\inn(Z)$. As $|\Aut(Z)/\inn(Z)|=2$,
by Lemma~\ref{facts about AGL23} (v), there are exactly two
$(A_X,A_Y)$-double cosets in $\Aut(Z)$. Now the Goldschmidt Lemma
implies that there are exactly two isomorphism classes of amalgams
of weak type $(X,Y,Z)$.
\end{proof}

By  Lemma~\ref{uniqueness}, if $\mathcal B$ is an amalgam of weak type $(X,Y,Z)$, then $\mathcal B$
is isomorphic to $\AA$ or $\AA'$. Let $F$ be the universal completion of $\AA$.

\begin{lemma}\label{presentation}
Let
\[\begin{array}{rcl}
 \mathcal{R}_Z&=&\{a^3=b^3=[a,b]^3=[a,[a,b]]=[b,[a,b]]=1, \\
        \;&\;& ~t^2=u^2=[t,u]=1,~a^u=a,~ b^u=b^{-1},~ a^t=a^{-1},~ b^t=b\};\\
\mathcal{R}_X&= & \{p^4=q^4=1,~p^2=q^2=t,~q^p= q^{-1},~
a^p=[a^{-1},b],~ [a,b]^q=[b,a]a,\\
\;&\;&\; p^u=p^{-1},~p^b=pq, ~q^b=p,~p^{b^2}=q\}; \mbox{ and}\\
\mathcal{R}_Y&=& \{r^4=s^4=1,~r^2=s^2=u,s^r = s^{-1}, ~
~[a^{-1},b]^r=b ,~[a,b]^s=b\inv [b,a],\\
\;&\;& \;r^t=r^{-1}, s^t = sr, ~r^a=rs,~ s^a=r,~r^{a^2}= s\}.\end{array}\]  Then $F\cong
\<a,b,p,q,r,s,t,u\mid \mathcal{R}_X,\mathcal{R}_Y,\mathcal{R}_Z \>$.
\end{lemma}
\begin{proof}

Let
$F^*=\<a,b,p,q,r,s,t,u\mid\mathcal{R}_X,\mathcal{R}_Y,\mathcal{R}_Z
\>$,  $Z^*= \<a,b,t,u\mid\mathcal{R}_Z\>$,
$X^*=\<a,b,p,q,t,u\mid\mathcal{R}_X,\mathcal{R}_Z \>$, and
$Y^*=\<a,b,r,s,t,u\mid\mathcal{R}_Y, \mathcal{R}_Z\>$. Define a map
$\Theta$ from $\{a,b,p,q,r,s,t,u\}$ into $\mathrm{SL}_3(3)$ as
follows:
\[ a \mapsto {\left( \begin{smallmatrix}
1 & 0 & 0 \\
1 & 1 & 0 \\
0 & 0 & 1 \end{smallmatrix} \right)},~~ b \mapsto {\left(
\begin{smallmatrix}
1 & 0 & 0 \\
0 & 1 & 0 \\
0 & 1 & 1 \end{smallmatrix} \right)},~~ t \mapsto {\left(
\begin{smallmatrix}
1 & 0 & 0 \\
0 & -1 & 0 \\
0 & 0 & -1 \end{smallmatrix} \right)},~~ u \mapsto {\left(
\begin{smallmatrix}
-1 & 0 & 0 \\
0 & -1 & 0 \\
0 & 0 & 1 \end{smallmatrix} \right)},\]

\[ p \mapsto {\left(\begin{smallmatrix}
1 & 0 & 0 \\
0 & 0 & 1 \\
0 & -1 & 0 \end{smallmatrix} \right)},~~ q \mapsto {\left(
\begin{smallmatrix}
1 & 0 & 0 \\
0 & -1 & 1 \\
0 & 1 & 1 \end{smallmatrix} \right)},~~
 r \mapsto {\left(\begin{smallmatrix}
0 & 1 & 0 \\
-1 & 0 & 0 \\
0 & 0 & 1 \end{smallmatrix} \right)},~~ s \mapsto {\left(
\begin{smallmatrix}
-1 & 1 & 0 \\
1 & 1 & 0 \\
0 & 0 & 1 \end{smallmatrix} \right).}\] Then
$\Theta_{|\{a,b,p,q,t,u\}}$ extends to a homomorphism from $X^*$ to
$X$, $\Theta_{|\{a,b,r,s,t,u\}}$ extends to a homomorphism from
$Y^*$ to $Y$ and $\Theta_{|\{a,b,t,u\}}$ extends to a homomorphism
from $Z^*$ to $Z$ (recall that $X$, $Y$ and $Z$ are subgroups of
$\SL_3(3)$). Therefore $X$ is a quotient of $X^*$, $Y$ is a quotient
of $Y^*$ and $Z$ is a quotient of $Z^*$. Coset enumeration using
{\sc Magma} \cite{magma}, for example,  gives the orders of $X^*$,
$Y^*$ and $Z^*$ and so we get $X\cong X^*$, $Y\cong Y^*$ and $Z\cong
Z^*$. The  obvious identifications $\phi_1:Z^*\rightarrow X^*$,
$\phi_2: Z^* \rightarrow Y^*$ are monomorphisms and so $\mathcal
A^*=(X^*,Y^*, Z^*,\phi_1,\phi_2)$ is an amalgam of weak type
$(X,Y,Z)$. Furthermore, the  identifications $\psi_1:X^*\rightarrow
F^*$, $\psi_2: Y^* \rightarrow F^*$ are monomorphisms and so $F^*$
is a completion of  $\A^*$. Moreover the map $\Theta$ defines
embeddings of $X^*$ and $Y^*$ in $\mathrm{SL}_3(3)$ and so it
follows that $\SL_3(3)$ is also a completion of $\A^*$. Lemma
\ref{A' has non-trivial O3} and Lemma \ref{uniqueness} therefore
imply that $\A^*$ is isomorphic to $\A$. Now using \cite[Theorem 1,
p.3]{serre} we deduce that $F^*$ is the universal completion of
$\A$. Hence $F \cong F^*$.
\end{proof}

Because of Lemma~\ref{presentation} we may and do now take
$F=\<a,b,p,q,r,s,t,u\mid \mathcal{R}_X,\mathcal{R}_Y,\mathcal{R}_Z
\>$,  $X=\<a,b,p,q,t,u\mid\mathcal{R}_X,\mathcal{R}_Z \>$,
$Y=\<a,b,r,s,t,u\mid\mathcal{R}_Y, \mathcal{R}_Z\>$ and $Z=
\<a,b,t,u\mid\mathcal{R}_Z\>$ identified as subgroups of $F$. We
mention that the presentation of $F$ in Lemma~\ref{presentation} is
by no means an efficient presentation but instead has been
intentionally constructed so that the subgroups
$O_2(C_{X}(t))=\<p,q\>$ and $O_2(C_{Y}(u))=\<r,s\>$, both isomorphic
to the quaternion group of order 8, are accessible. Notice that the
involutions $t$ and $u$ satisfy $ O_3(X) t\in Z(X/O_3(X))$ and
$O_3(Y)u \in Z(Y/O_3(Y))$.

\begin{lemma}\label{2Q8's} Let $P^* = \<p,q\>$, $R^* = \<r,s\>^{pr}$. Then $P^*\cong R^*\cong Q_8$,
$\langle b,u\rangle \cong \sym(3)$ normalizes both $P^*$ and $R^*$ and  $\langle
P^*,R^*\>\<b,u\rangle \le C_F(t)$. Moreover, $P^*\langle b,u\rangle  \cong R^*\langle b,u\rangle
\cong \GL_2(3)$.
\end{lemma}
\begin{proof}
We have already said that $\<p,q\>$ and $\langle r,s\>$ are
isomorphic to $\Q_8$. That $\<b,u\>\cong\sym(3)$ and normalizes
$P^*$ is evident from the relations in $\mathcal R_Z \cup \mathcal
R_X$. Furthermore,  $P^*\<b,u\> = C_X(t)\cong \GL_2(3)$ follows
directly from the structure of $X$. Now  using the relations $p^u=
p^{-1}$ and $r^t= r^{-1}$, we get
$$u^{pr} = r^{-1}p^{-1}upr = r^{-1}u(p^{-1})^{u}pr = r^{-1}up^2r = r^{-1}utr=r^{-1}ur^t t = t.$$
Thus $t \in Z(R^*)$ and so $R^* \le C_F(t)$. We now observe that
$a^{pr} = [a^{-1},b]^r = b$ and
$(ut)^{pr}=u^{pr}t^{pr}=tt^{r}=tut=u$ follow from the relations in
$\mathcal R_X \cup \mathcal R_Y \cup \mathcal R_Z$. Therefore
$\<R^*,b,u\>=\<r,s,a,ut\>^{pr} \leq Y^{pr}$. Since $\<a,ut\>$
normalizes $\<r,s\>$, we see that $\<b,u\>$ normalizes $R^*$.
Finally, a calculation in $Y$ proves that $\mathrm{GL}_2(3) \cong
\<r,s,a,ut\>\cong \<R^*,b,u\>$.
\end{proof}

\section{Proof of  Theorem A}

We now work under the hypothesis of Theorem~\ref{maintheorem}A. So
we assume that $G$ is a finite group generated by distinct subgroups
$A$ and $B$ with  $A \cong B\cong \AGL_2(3)$. Further we assume that
$N_G(Z(S))\leq A$ where $S \in \mathrm{Syl}_3(A\cap B)$. We continue
the notation from the last section. In particular, we recall that
$F$ is the universal completion of the amalgam $\AA$. Set $C = A\cap
B$ and let $S \in \syl_3(C)$.

\begin{lemma}\label{lem1} The following hold:
\begin{enumerate}[$(i)$]
\item $S \in \syl_3(A) \cap \syl_3(B) \subseteq \syl_3(G)$;
\item $C=N_A(S)=N_B(S)=N_G(S)$; and
\item  $G$ is a completion of $\A$.\end{enumerate}
\end{lemma}
\begin{proof} Since $S \le C \le B$, there exists a Sylow $3$-subgroup $S_1$ of $B$ which contains
$S$. Now $Z(S)$ is a characteristic subgroup of $S$ and so is normal in $N_{S_1}(S)$. As $N_G(Z(S))
\le A$, it follows that $N_{S_1}(S)\le A \cap B = C$. Because $S \in \syl_3(C)$, we deduce that
$S_1 = S \in \syl_3(B)$. As $A \cong \B$, we  have $S \in \syl_3(A)$. Since $N_G(S)\le N_G(Z(S))=
N_A(S)$, we also have $S \in \syl_3(G)$. Thus (i) holds.

Now $N_B(S) \le N_G(Z(S)) = N_A(Z(S))= N_A(S)$ by Lemma~\ref{facts
about AGL23} (iv) and so $N_A(S)=N_B(S)\leq C$. Since $A \neq B$ and
$N_A(S) $ is a maximal subgroup of $A$ by Lemma~\ref{facts about
AGL23} (iv) again, we have $C=N_A(S)=N_B(S)$.

Assume that $O_3(G) \neq 1$. Then $O_3(G) \le S$ by (i). In
particular, $O_3(G)$ is normal in $A$. Hence $O_3(G) = O_3(A)$ has
order $9$ by Lemma \ref{facts about AGL23} (i) and so contains
$Z(S)$. Since $C_G(O_3(G)) \le N_G(Z(S)) \le A$, we have
$C_G(O_3(G)) = O_3(G)$ by Lemma~\ref{facts about AGL23} (i). It
follows that $G/O_3(G)$ is isomorphic to a subgroup of $\GL_2(3)$
and so $G$ has order at most 432. But then  $A= B$, which is a
contradiction. Thus $O_3(G)=1$. So we have that $G$ is a completion
of an amalgam of weak type $(X,Y,Z)$ and it therefore follows from
Lemma~\ref{A' has non-trivial O3}  and Lemma~\ref{uniqueness}  that
$G$ is a completion of $\A$.
\end{proof}

By Lemma~\ref{lem1}, there is a homomorphism $\Psi:F \rightarrow G$
such that $X\Psi = A$, $Y\Psi= B$ and $Z\Psi = C$. We will denote
the images  of elements of $F$ under $\Psi$  by boldface letters.
Thus for example $a \Psi= \mathbf a$. Define $P= P^*\Psi= \<{\mathbf
p},{\mathbf q}\>$, $R = R^*\Psi= \<{\mathbf r} ,\mathbf s
\>^{{\mathbf p}{\mathbf r}}$ and $K = \<\mathbf b,\mathbf u\>$.
Notice that, as $\Psi$ restricted to both $X$ and $Y$ is a
monomorphism, we have $PK \cong RK \cong \GL_2(3)$ by
Lemma~\ref{2Q8's}. Set $Q_A=O_3(A)$ and $Q_B=O_3(B)$. Then $\mathbf
t Q_A \in Z(A/Q_A)$. Let $H = C_G(\mathbf t)$ and $\overline H =
H/\<{\mathbf t}\>$.

\def \ov {\overline}
\begin{lemma}\label{involution centralizer-SC3}We have  $C_{H}({\mathbf b} )= \<{\mathbf b} ,{\mathbf t}\>$.
\end{lemma}
\begin{proof}
We see directly from the presentation of $F$, that $b \in O_3(Y)$.
It follows that ${\mathbf b} \in Q_B$ and consequently, using
Lemma~\ref{facts about AGL23} (i), ${\mathbf b} \in Z(T)$ for some
Sylow $3$-subgroup $T$ of $B$. Since $T \in\syl_3(G)$ by
Lemma~\ref{lem1} (i),  $C_G({\mathbf b}) \le N_G(\<{\mathbf b}\>)
\le B$. Notice also that $t \in Y$ and so ${\mathbf t} \in B$. We
now calculate that directly in $Y$ to see that  $C_G(\<{\mathbf
b},{\mathbf t}\>) = C_B(\<{\mathbf b},{\mathbf t}\>)=\<{\mathbf
b},{\mathbf t}\>$.
\end{proof}
Set $W=\<P,R\>$. By Lemma~\ref{2Q8's}, $W \le H$. Also, as $K$ normalizes $P$ and $R$, $K$
normalizes $W$.

\begin{lemma}\label{lastlem} We have that $W$ is a
$2$-group which is normalized by $K$. Furthermore, $C_W({\mathbf
b})= \<{\mathbf t}\>$.
\end{lemma}
\begin{proof}

 Lemma~\ref{involution centralizer-SC3} implies that we may apply Theorem~\ref{Feit-Thompson} to $\ov H$. If
$\ov H \cong \mathrm{PSL}_2(7)$, then as $\Q_8 \cong R \le H$, $H
\cong \mathrm{SL}_2(7)$. But then ${\mathbf t}$ is the unique
involution in $H$ contrary to  $\mathbf u \neq {\mathbf t}$ and
$\mathbf u \in H$. So $\ov H$ contains a nilpotent normal subgroup
$\ov {N}$ such that $\ov{H}/\ov{N}\cong \mathrm{Alt}(5)$ or
$\sym(3)$ or $C_3$. Therefore $H$ contains a nilpotent normal
subgroup $3'$-subgroup $N$ containing ${\mathbf t}$ such that
$C/N\cong \mathrm{Alt}(5)$ or $\sym(3)$ or $C_3$. Now we claim that
$P$ and $R$ are both contained in $N$. By Lemma~\ref{2Q8's},
$PK\cong \mathrm{GL}_2(3)$ has two normal $3'$-subgroups containing
$\mathbf{t}$. One is $\<\mathbf{t}\>$ the other is $P\cong Q_8$. If
$PK \cap N$ has order 2 then $PK/(PK\cap N)\cong PKN/N$ has order a
multiple of 8 which is not possible. Therefore $PK \cap N=P$. By the
same argument $R \le N$ and since $N$ is nilpotent, $W$ is a
$2$-subgroup of $N$. That $C_W({\mathbf b})= \<{\mathbf t}\>$
follows from Lemma~\ref{involution centralizer-SC3}.
\end{proof}

Define $F_1 = \<F\mid\mathcal{R}_1\>$, $F_2 =
\<F\mid\mathcal{R}_2\>$, $F_3 = \<F\mid\mathcal{R}_3\>$ and $F_4 =
\<F\mid \mathcal{R}_4\>$ where, for $x = r^{pr}$ and $y = s^{pr}$,
$\mathcal{R}_1=\{[x,q][p,y]=1, [q,x][xy,pq]=1\}$,
$\mathcal{R}_2=\{[x,q][p,y]=t, [q,x][xy,pq]=1\}$,
$\mathcal{R}_3=\{[x,q][p,y]=1, [q,x][xy,pq]=t\}$ and
$\mathcal{R}_4=\{[x,q][p,y]=t, [q,x][xy,pq]=t\}$. The next lemma
completes the proof of Theorem~A.

\begin{lemma}
$G \cong \mathrm{M}_{12}$ or $\mathrm{SL}_3(3)$.
\end{lemma}
\begin{proof} We first calculate in $F$. Set $x = r^{pr}$ and $y = s^{pr}$.  The relation $p^{b^2}=
q$ is in $\mathcal R_X$ and,  using the equality $b = a^{pr}$ used
in Lemma~\ref{2Q8's}, we have $x^{b^2}= r^{prb^2}=
r^{aprb}=r^{a^2pr}=(rs)^{apr}=(rsr)^{pr}=s^{pr}=y.$ Thus in $F$, we
have the equality
$$[x,q][p,y]=[x,p^{b^2}][p,x^{b^2}].$$ Similarly, using the relations in $\mathcal R_X \cup R_Y$, we
have $ (pq)^b = q$ and $x^b = xy$ and so
$$[q,x][xy,pq]=[(pq)^b,x][x^b,pq].$$ Applying $\Psi$, we now have
$[{\mathbf x},{\mathbf q}][{\mathbf p},{\mathbf y}]=[{\mathbf
x},{\mathbf p}^{{\mathbf b}^2}][{\mathbf p},{\mathbf x}^{{\mathbf
b}^2}]$ and $[{\mathbf q},{\mathbf x}][{\mathbf x}{\mathbf
y},{\mathbf p}{\mathbf q}]=[({\mathbf p}{\mathbf q})^{\mathbf
b},{\mathbf x}][{\mathbf x}^{\mathbf b},{\mathbf p}{\mathbf q}].$
Since, by Lemma~\ref{lastlem}, $W=\<{\mathbf p},{\mathbf q},{\mathbf
x},{\mathbf y}\>$ is a $2$-group and $C_W({\mathbf b}) =\<{\mathbf
t}\>$, we may apply Burnside's Lemma to $W/\<{\mathbf t}\>$ (with
$\zeta$ as conjugation by ${\mathbf b}^2$), to get the following
equality of cosets $\<{\mathbf t}\>[{\mathbf x},{\mathbf
p}^{{\mathbf b}^2}] =\<{\mathbf t}\>[{\mathbf x},{\mathbf
p}]^{\mathbf b}$ and $\<{\mathbf t}\>[{\mathbf p},{\mathbf
x}^{{\mathbf b}^2}]=\<{\mathbf t}\>[{\mathbf p},{\mathbf
x}]^{\mathbf b}.$ Therefore
$$\<{\mathbf t}\>[{\mathbf x},{\mathbf q}][{\mathbf p},{\mathbf y}]=
\<{\mathbf t}\>[{\mathbf x},{\mathbf p}^{{\mathbf b}^2}][{\mathbf
p},{\mathbf x}^{{\mathbf b}^2}]= \<{\mathbf t}\>[{\mathbf
x},{\mathbf p}]^{\mathbf b}[{\mathbf p},{\mathbf x}]^{\mathbf b} =
\<{\mathbf t}\>.$$ Hence $ [{\mathbf x},{\mathbf q}][{\mathbf
p},{\mathbf y}] \in \<{\mathbf t}\>$. Similarly,$$ \<{\mathbf
t}\>[{\mathbf q},{\mathbf x}][{\mathbf x}{\mathbf y},{\mathbf
p}{\mathbf q}]=\<{\mathbf t}\>[({\mathbf p}{\mathbf q})^{\mathbf
b},{\mathbf x}][{\mathbf x}^{\mathbf b},{\mathbf p}{\mathbf q}]=
\<{\mathbf t}\>[{\mathbf p}{\mathbf q},{\mathbf x}]^{{\mathbf
b}^2}[{\mathbf x},{\mathbf p}{\mathbf q}]^{{\mathbf b}^2}=\<{\mathbf
t}\>$$ and so
 $[{\mathbf q},{\mathbf x}][{\mathbf x}{\mathbf y},{\mathbf p}{\mathbf q}]\in \<{\mathbf t}\>$. It follows that
$G$ is a quotient of one of the groups $F_1$, $F_2$, $F_3$ or $F_4$.
Coset enumeration, using for example {\sc Magma} \cite{magma}, gives
$|F_1|=95040$, $|F_2|= 1$, $|F_3|= 5616$ and $|F_4|= 1.$ Since both
$\mathrm{M}_{12}$ and $\mathrm{SL}_3(3)$ satisfy the hypothesis of
Theorem~\ref{maintheorem}A, they are both completions of $\mathcal
A$ and hence are quotients of one of $F_1$ or  $F_3$. By comparing
the orders of the groups we infer that $F_1 \cong \mathrm{M}_{12}$
and $F_3 \cong \mathrm{SL}_3(3)$. Hence $G \cong \M_{12}$ or $G\cong
\SL_3(3)$ as claimed.
\end{proof}

\section{Proof of Theorem B }\label{graphsection}

Suppose that $\Gamma=(V,E)$ is a graph and $G$ is a group of automorphisms of $\Gamma$. Then, for
$\a \in V$, $\Gamma(\a)$ is the set of neighbours of $\a$,  $G_\a = \stab_G(\a)$ and, for
$\{\a,\b\}\in E$, $G_{\a\b}=\stab_G(\{\a,\b\})$. For  $H \le G$, we denote the subgraph of $\Gamma$
fixed by $H$ by $\Gamma^H$.

Now suppose that $G$ is a group and $X$ and $Y$ are subgroups of $G$. Define $\Gamma=\Gamma(G,X,Y)$
to be the bipartite graph with vertex set $V(\Gamma)=\{Xg|g \in G\} \cup\{Yg\mid g\in G\}$ and edge
set $E(\Gamma)=\{\{Xg,Yh\}\mid Xg \cap Yh \neq \emptyset\}$. Then $\Gamma$ is a graph which admits
an action of $G$ which is edge but not vertex transitive.

\begin{lemma} \label{pregraph}Suppose that $\Gamma=(V,E)$ is a  graph with no vertex of degree one and $G$ is a group of automorphisms of
$\Gamma$ which acts transitively on $E$ but not on $V$. Then, for $\{\a,\b\}\in E$ and $\theta \in
\{\a,\b\}$,
\begin{enumerate}[$(i)$]
\item $G$ has exactly two orbits $\a G$ and $\b G$ on $V$;
\item $G_\theta$ acts  transitively on $\Gamma(\theta)$;
\item $\Gamma \cong \Gamma(G,G_\a,G_\b)$; and
\item $G=\<G_\a,G_\b\>$ if and only if $\Gamma$ is connected.
\end{enumerate}
\end{lemma}

\begin{proof}  (i) Since $G$ acts transitively on $E$ and every vertex is contained in an edge, $G$ has at most two orbits
on $V$. As $G$ is not transitive on $V$, we deduce that the orbits of $G$ on $V$ are $\a G$ and $\b
G$.

(ii) Let $\mu, \tau \in \Gamma(\theta)$. Then, as $G$ is transitive on $E$, there is an element of
$G$ moving $\{\theta,\mu\}$ to $\{\theta,\tau\}$. By (i), this element fixes $\theta$. Hence
$G_\theta$ acts transitively on $\Gamma(\theta)$.

(iii) We define a map $ \Gamma(G,G_\a, G_\b) \rightarrow \Gamma$ by $G_\theta h \mapsto \theta h$
for $\theta \in \{\a,\b\}$ and $h \in G$. This map is well-defined on both the vertices and edges
of $\G(G,G_\a,G_\b)$ and is easily checked to be a bijection.

(iv) This is \cite[2.4]{goldschmidt}.

\end{proof}

Recall the graphs $\Gamma_1$ and $\Gamma_2$ described in the introduction. Then, by
Lemma~\ref{pregraph} (iii),  $\Gamma_1\cong \Gamma(\mathrm{PSL}_3(3),A_1,B_1)$ and $\Gamma_2\cong
\Gamma(\mathrm{M}_{12},A_2,B_2)$ .

We now prove Theorem~B.

\begin{proof} [Proof of Theorem~B]Assume that $G$ and $\Gamma=(V,E)$ are as in the statement of Theorem~B. Let $z$ be a non-trivial element of $O_3(G_\a)$.
Then, by hypothesis, $\Gamma^z$ is a tree with at least one edge. As
$\a\in \Gamma^z$, we may pick an edge $\{\a,\b\} \in \Gamma^z$. Let
$T_0 \in \syl_3(G_{\a\b})$ with $z\in T_0$ and set $T = Z(T_0)$.
Notice that either $z \in T$ or $T_0$ is non abelian and then
$T_0\in \syl_3(G_\a)$ follows from Lemma \ref{facts about AGL23}
(iii). In either case $T \cap O_3(G_\a)\neq 1$. Let $z'$ be a
non-trivial element of $T \cap O_3(G_\a)$. Then $\Gamma^T\subseteq
\Gamma^{z'}$ which is a tree. Hence $\Gamma^T$ is also a tree and
$\Gamma^T$ contains  $\{\a,\b\}$. Let
 $H=N_G(T)$. Then $H$ acts on $\Gamma^T$. Since $H$ is finite and $\Gamma^T$ is a
tree, $H$ fixes either a vertex or an edge of $\Gamma^T$. Using
Lemma~\ref{pregraph}(i), we get that in either case $H$ fixes a
vertex $\theta$ in $\Gamma^T$. Therefore $N_{G_ \a}(T) \le H
=N_{G_\theta}(T)$. If $|T|=3$, then, as $G_\theta$ has extraspecial
Sylow $3$-subgroups of order $27$, either $T=T_0=\<z\>$ or
$|T_0|=27$. Hence in either case we have that $T \le O_3(G_\a)$ and
$N_{G_\a}(T)$ has index $4$ in $G_\a$. Therefore $H$ has index at
most $4$ in $G_\theta$ and, as $H$ normalizes $T$ which has order
$3$, we deduce that $H = N_{G_\a}(T)$. Suppose that $|T|>3$, then
$T_0= T$ has order $9$. It follows that $N_{G_\a}(T)$ contains a
Sylow $3$-subgroup $T_1$ of $G_\alpha$. Since $T_1\le N_{G_\a}(T)
\le H \le G_{\theta}$, we see that $T_1$ fixes the unique path in
$\Gamma^T$ from $\a$ to $\theta$. If $\alpha\neq \theta$, then $T_1$
fixes an edge in $\Gamma^T$. However $G$ is transitive on edges and
so each edge stabilizer is conjugate in $G$. This contradicts
$T=T_0$ being a Sylow $3$-subgroup of $G_{\a\b}$. Therefore $\alpha
= \theta$. In conclusion we have $H \le G_\a$. Set $A = G_\a$, $B =
G_\b$ and $S = T_0$. Then, as $\Gamma$ is connected,  $G=
\<G_\a,G_\b\rangle=\<A,B\>$ by Lemma~\ref{pregraph} (iv).
Furthermore $A \neq B$ by Lemma~\ref{pregraph} (ii), $A\cong B\cong
\AGL_2(3)$ and $N_G(Z(S))=H\le A$. Thus Theorem~A is applicable and
gives  $G \cong \M_{12}$ or $G\cong \SL_3(3)$. It now follows that
$\Gamma\cong \Gamma_1= \Gamma(\SL_3(3),A_1,B_1)$ or $\Gamma\cong
\Gamma_2=\Gamma(\M_{12},A_2,B_2)$.
\end{proof}

\bibliographystyle{plain}
\bibliography{mybibliography}
\end{document}